\documentclass[11pt,epsf]{article}
\usepackage{amsthm}
\usepackage{amsfonts}
\usepackage{amssymb}
\title{A property of a partial theta function}
\author{Vladimir Petrov Kostov\\ 
Universit\'e de Nice, 
Laboratoire de Math\'ematiques, Parc Valrose,\\ 06108 Nice Cedex 2, France,  
e-mail: kostov@math.unice.fr} 
\date{}
\newtheorem{tm}{Theorem}

\newtheorem{rem}[tm]{Remark}
\newtheorem{rems}[tm]{Remarks}
\newtheorem{lm}[tm]{Lemma}

\begin{document} 
\maketitle 

\begin{abstract}
The series $\theta (q,x):=\sum _{j=0}^{\infty}q^{j(j+1)/2}x^j$ 
converges for $|q|<1$  
and defines a {\em partial theta function}. For any fixed 
$q\in (0,1)$ it has infinitely many negative zeros. 
It is known that for $q$ taking one of the {\em spectral} values $\tilde{q}_1$, 
$\tilde{q}_2$, $\ldots$   
(where $0.3092493386\ldots =\tilde{q}_1<\tilde{q}_2<\cdots <1$, 
$\lim _{j\rightarrow \infty}\tilde{q}_j=1$) the function $\theta (q,.)$ 
has a double zero which is the rightmost of its real zeros 
(the rest of them being simple). 
For $q\neq \tilde{q}_j$ the partial theta function has no multiple real zeros. 
We prove that: 

1) for $q\in (\tilde{q}_{j},\tilde{q}_{j+1}]$ the function 
$\theta$ is a product of a degree $2j$ real polynomial without real roots 
and a function of the Laguerre-P\'olya class $\cal{LP-I}$;

2) for $q\in \mathbb{C}\backslash 0$, $|q|<1$,    
$\theta (q,x)=\prod _i(1+x/x_i)$, where $-x_i$ are the zeros of $\theta$; 
%counted with multiplicity;

3) for any fixed $q\in \mathbb{C}\backslash 0$, $|q|<1$, the function $\theta$ has 
at most finitely-many multiple zeros;

4) for any $q\in (-1,0)$ the function $\theta$ is a product of a 
real polynomial without real zeros and a 
function of the Laguerre-P\'olya class $\cal{LP}$. 

5) for any fixed $q\in \mathbb{C}\backslash 0$, $|q|<1$, and for $k$ sufficiently 
large, the function $\theta$ has a zero $\zeta _k$ close to $-q^{-k}$.   
These are all but finitely-many of the zeros of $\theta$. 
%The same is true for $q\in \mathbb{C}\backslash 0$ with $|q|$ 
%sufficiently small and for any $k$.
\\ 


\end{abstract}

\section{Introduction}

The series in two variables $\theta (q,x):=\sum _{j=0}^{\infty}q^{j(j+1)/2}x^j$ 
defines an entire function in $x$ for each fixed $q$, $|q|<1$. 
In the present paper we consider three situations: 
$q\in \mathbb{C}$, $|q|<1$, $x\in \mathbb{C}$ or 
$q\in [0,1)$, $x\in \mathbb{R}$ or $q\in (-1,0]$, $x\in \mathbb{R}$. 
We say that the series defines a 
{\em partial theta function} (the Jacobi theta function is defined 
by the series $\Theta (q,x):=\sum _{j=-\infty}^{\infty}q^{j^2}x^j$ and one has 
$\theta (q^2,x/q)=\sum _{j=0}^{\infty}q^{j^2}x^j$, i.e. in the definition of 
$\theta$ only partial summation is performed). We regard $q$ as a 
parameter and $x$ as a variable. 

There are several domains in which the partial theta function is an object of 
interest: the theory 
of (mock) modular forms (\cite{BrFoRh}), asymptotic analysis (\cite{BeKi}),
statistical physics 
and combinatorics (\cite{So}) and Ramanujan type $q$-series 
(\cite{Wa}). Additional information about $\theta$ can be found in \cite{AnBe}. 

A recent study relates the function $\theta$ to a problem considered by 
 Hardy, Petrovitch and Hutchinson, see \cite{Ha}, \cite{Hu}, \cite{Pe}, 
\cite{KaLoVi}, \cite{KoSh} 
and \cite{Ko2}. In the article \cite{KaLoVi} the existence of a constant 
$\tilde{q}\in (0,1)$ is proved such that for $q\in (0,\tilde{q})$ 
the function $\theta (q,.)$ 
has only real negative simple zeros while $\theta (\tilde{q},.)$ has a 
double negative zero the rest of the zeros being negative and simple. The more 
accurate value $0.3092493386\ldots$ of $\tilde{q}$ is given in \cite{KoSh}. 
The function $\theta$ belongs to the Laguerre-P\'olya class $\cal{LP-I}$ 
exactly for $q\in (0,\tilde{q}]$. We remind that a function of this class 
is either representable in the form 

\begin{equation}
\psi (x) = c x^{m} e^{\sigma x} \prod_{k=1}^{\omega} (1+x/x_{k}),
\label{LPI}
\end{equation} 
where $\omega \in \mathbb{N}\cup \infty$, $c\in \mathbb{R}$, 
$\sigma \geq 0$, $x_k>0$, $m\in \mathbb{N}\cup 0$ and $\sum 1/x_{k} < \infty$ 
or $\psi (-x)$ is representable in this form. 
For each compact set in $\mathbb{C}$, the restriction to it of 
each function of this class is a uniform 
limit of polynomials with all roots real and negative (or real and positive). 

An entire function belongs 
to the Laguerre-P\'olya class $\cal{LP}$ if it is of the form 

\begin{equation}
\psi (x) = c x^{m} e^{-\alpha x^{2} + \beta x} \prod_{k=1}^{\omega} (1+x/x_{k})e^{-x/x_{k}},
\label{LP}
\end{equation}
where $\omega \in \mathbb{N}\cup \infty$, $c, \beta \in \mathbb{R}$, 
$\alpha \geq 0$, $m\in \mathbb{N}\cup 0$ and $\sum x_{k}^{-2} < \infty$.
For each compact set in $\mathbb{C}$, the restriction to it of 
each function of this class is a uniform 
limit of polynomials with all roots real. 

The {\em spectrum} of $\theta$ is the set $\Gamma$ 
of values of $q$ for which $\theta (q,.)$ has a multiple real zero. 
(This terminology has been proposed by B.Z.~Shapiro, see \cite{KoSh}.)  
It is shown in \cite{Ko2} that $\Gamma$ consists of countably-many real 
numbers denoted by 
$0<\tilde{q}=\tilde{q}_1<\tilde{q}_2<\cdots <\tilde{q}_j<\cdots <1$, 
$\lim _{j\rightarrow \infty}\tilde{q}_j=1$. 
For $\tilde{q}_j\in \Gamma$ the function $\theta (\tilde{q}_j,.)$ 
has exactly one multiple  
real zero which is negative, of multiplicity $2$ and is the rightmost 
of its real zeros.

In the present paper we prove the following theorem:

\begin{tm}\label{basictm}
(1) For any fixed $q\in \mathbb{C}\backslash 0$, $|q|<1$, and for $k$ sufficiently 
large, the function $\theta (q,.)$ has a zero $\zeta _k$ close to $-q^{-k}$ 
(in the sense that $|\zeta _k+q^{-k}|\rightarrow 0$ as $k\rightarrow \infty$). 
These are all but finitely-many of the zeros of $\theta$.

%(2) The same statement holds true for $q\in \mathbb{C}\backslash 0$ with $|q|$ 
%sufficiently small and for any $k\in \mathbb{N}$. 

(2) For any $q\in \mathbb{C}\backslash 0$, $|q|<1$, one has 
$\theta (q,x)=\prod _k(1+x/x_k)$, where $-x_k$ are the zeros of $\theta$ 
counted with multiplicity.

(3) For $q\in (\tilde{q}_{j},\tilde{q}_{j+1}]$ the function 
$\theta (q,.)$ is a product of a degree $2j$ real polynomial without real roots 
and a function of the Laguerre-P\'olya class $\cal{LP-I}$. 
Their respective forms are $\prod _{k=1}^{2j}(1+x/\eta _k)$ and 
$\prod _k(1+x/\xi _k)$, where $-\eta _k$ and $-\xi _k$ are  
the complex and the real zeros of $\theta$ counted with multiplicity.

(4) For any fixed $q\in \mathbb{C}\backslash 0$, $|q|<1$, the function $\theta (q,.)$ has 
at most finitely-many multiple zeros.

(5) For any $q\in (-1,0)$ the function $\theta (q,.)$ is a product of the form 
$R(q,.)\Lambda (q,.)$, where $R=\prod _{k=1}^{2j}(1+x/\tilde{\eta}_k)$ 
is a real polynomial with constant term $1$ and without real zeros and 
$\Lambda =\prod _k(1+x/\tilde{\xi}_k)$, $\tilde{\xi}_k\in \mathbb{R}^*$, is a 
function of the Laguerre-P\'olya class $\cal{LP}$. One has 
$\tilde{\xi}_k\tilde{\xi}_{k+1}<0$. The sequence $\{ |\tilde{\xi}_k|\}$ is 
monotone increasing for $k$ large enough.
\end{tm}

\begin{rems}
{\rm (1) In the present paper differentiation w.r.t. $q$ is not used, 
therefore for a function $f(q,x)$ we write $f'(q,x)$ 
instead of $(\partial f/\partial x)(q,x)$.

(2) To prove part (1) of the theorem 
we use properties of the (true) Jacobi theta function 
$\Theta (q,x):=\sum _{j=-\infty}^{\infty}q^{j^2}x^j$. It is known that $\Theta$ has 
only simple zeros (see \cite{Wi}). Hence this is also the case of the function 
$\Theta ^*(q,x)=\Theta (\sqrt{q},\sqrt{q}x)=
\sum _{j=-\infty}^{\infty}q^{j(j+1)/2}x^j$. 
The zeros of $\Theta ^*$ are all the numbers $\mu _k:=-q^{-k}$, 
$k\in \mathbb{Z}$ 
(this follows from the formula for the zeros of $\Theta$, see \cite{Wi}). 
One has 

\begin{equation}\label{Theta}
\Theta ^*(q,x)=qx\Theta ^*(q,qx)~.
\end{equation} 
The function $\Theta$ satisfies the following identity 
(known as the Jacobi triple product, see~\cite{Wi}): 

%\begin{equation}\label{Jacobi}
$$\Theta (q,x^2)=\prod _{m=1}^{\infty}(1-q^{2m})(1+x^2q^{2m-1})(1+x^{-2}q^{2m-1})$$
%\end{equation}
This implies the equality

\begin{equation}\label{Jacobi}
\Theta ^*(q,x^2)=\prod _{m=1}^{\infty}(1-q^m)(1+x^2q^m)(1+x^{-2}q^{m-1})
\end{equation}

(3) We write $x\in \Omega _k(\delta )$ instead of $|x-\mu _k|\leq \delta$, 
$\delta >0$.} 
%Hence $x\in \Omega _k(\delta )$ if and only if 
%$x/q\in \Omega _{k+1}(\delta /|q|)$.}
\end{rems}

\section{Proof of Theorem~\protect\ref{basictm}} 

%\begin{proof}[Proof of Theorem~\protect\ref{basictm}]

Represent the function $\Theta ^*$ in a neighbourhood of its zero 
$\mu _k=-q^{-k}$ in the form

$$\Theta ^*(q,x)=K(x+q^{-k})\tau (q,x)~~,~~{\rm where}~~
\tau =1+d_1(x+q^{-k})+d_2(x+q^{-k})^2+\cdots ~~,~~
K\in \mathbb{C}\backslash 0~~,~~d_j\in \mathbb{C}~.$$
Noticing that $x=x+q^{-k-s}-q^{-k-s}=-q^{-k-s}(1-q^{k+s}(x+q^{-k-s}))$ and 
using $s$ times equation (\ref{Theta}) 
one finds the following representation of 
$\Theta ^*$ in a neighbourhood of $\mu _{k+s}$:

$$\begin{array}{rcl}\Theta ^*(q,x)&=&q^{s(s+1)/2}x^sKq^s(x+q^{-k-s})
\tau (q,q^sx)\\ \\ 
&=&
K(-1)^sq^{s(-2k-s+3)/2}(1-q^{k+s}(x+q^{-k-s}))^s(x+q^{-k-s})
\tau (q,q^sx)~,~~{\rm where}\\ \\ 
\tau (q,q^sx)&=&1+d_1q^s(x+q^{-k-s})+d_2q^{2s}
(x+q^{-k-s})^2+\cdots ~.\end{array}$$
Introduce the new variable $X:=x+q^{-k-s}$. 
If one assumes that for all $s\in \mathbb{N}\cup 0$ the variable $X$ takes 
values in one and the same disk $|X|\leq \delta$, then as $s\rightarrow \infty$ 
the functions $(1-q^{k+s}X)^s$ and $\tau (q,q^sx)$ tend uniformly 
to $1$. Thus for $|X|\leq \delta$ (i.e. for $x\in \Omega _{k+s}(\delta )$) 
one has 

$$\Theta ^*=K(-1)^sq^{s(-2k-s+3)/2}(1+o(1))X~.$$ 
For any $B>0$ there exists $k_0\in \mathbb{N}\cup 0$ such that 
for $k\geq k_0$ one has $|q^{s(-2k-s+3)/2}|\geq B$ for any $s\in \mathbb{N}$. 
%In what follows we set $k=k_0$ (see the lemma). 
Hence there exists $\kappa >0$ such that for $|\eta |\leq \kappa$, 
for $k\geq k_0$ and for any 
$s\in \mathbb{N}$ the equation $\Theta ^*=\eta$ has a unique 
solution $X=X(\eta )$ with $|X|\leq \delta$.

Set $\theta (q,x)=\Theta ^*(q,x)+\Xi (q,x)$, where 
$\Xi (q,x)=-\sum _{j=-\infty}^{-1}q^{j(j+1)/2}x^j$. For fixed $q$, $|q|<1$, 
the series of $\Xi$ and $\Xi '$ 
both converge for $|x|>1$, and for any $\varepsilon >0$ 
there exists $A\geq 1$ such that if $|x|\geq A$, 
then $|\Xi|\leq \varepsilon$ and 
$|\Xi '|\leq \varepsilon$. (Indeed, 
both series as series of $1/x$ are without constant term and the moduli of 
all coefficients of the series of $\Xi$ are less than $1$.)

For $k\in \mathbb{N}$ sufficiently large the equation 

$$\theta (q,x)=0~~~~{\rm ,~~i.e.}~~~~\Theta ^*(q,x)=-\Xi (q,x)$$
has a unique solution $x=x(q)\in \Omega _k(\delta )$. Indeed, for such $k$ 
and for $x\in \Omega _k(\delta )$ 
%by the implicit function theorem, the above lemma and Theorem~\ref{tmGo}
the equation $\Theta ^*(q,x)=\eta$ 
has a solution $x(q,\eta )=\mu _k+O(\eta )$ which is holomorphic in $\eta$  
for $|\eta |$ sufficiently small. Substituting $-\Xi (q,x)$ for $\eta$ 
one obtains an equation of the form 
$x=\mu _k+\Delta (q,x)$, where $\max _{x\in \Omega _k(\delta)}|\Delta (q,x)|$ and 
$\max _{x\in \Omega _k(\delta)}|\Delta '(q,x)|$ 
can be made arbitrarily small by choosing $A$ and $k_0$ sufficiently large. 
Hence 
%(again by the implicit function theorem and Theorem~\ref{tmGo}) 
this equation has a unique solution in $\Omega _k(\delta)$.
%can be solved w.r.t. $x$. Lemma~\ref{greatder} implies that 
%this solution is unique.  

To complete the proof of part (1) we show that for $k_0$ large enough 
there remain only finitely-many zeros of $\theta$ outside the set 
$\Sigma :=\cup _{k=k_0}^{\infty}\Omega _k(\delta )$.  
We use equation 
(\ref{Jacobi}). (From now till the end of the proof of part (1) 
we consider $\theta (q,x^2)$ instead of $\theta (q,x)$, similarly for 
$\Theta$, $\Theta ^*$ and $\Xi$ which is not restrictive.) 
Suppose that $|x|^2\geq |q|^{-6}$ 
(this is not a restriction since 
the function $\theta$ is holomorphic and hence has 
finitely-many zeros in any compact set). Then for the products 
$\Pi _1:=\prod _{m=1}^{\infty}(1-q^m)$ and 
$\Pi _2:=\prod _{m=1}^{\infty}(1+x^{-2}q^{m-1})$ one has 

$$|\Pi _1|\geq \prod _{m=1}^{\infty}(1-|q|^m)=:r>0~~,~~
|x^{-2}q^{m-1}|<|q^m|~~{\rm and}~~ 
|\Pi _2|\geq \prod _{m=1}^{\infty}(1-|q|^m)=r~.$$ 
%$$|\Pi _2|=|1+x^{-2}|\prod _{m=2}^{\infty}|1+x^{-2}q^{m-1}|\geq (3/4)
%\prod _{m=1}^{\infty}(1-|q|^m)=3r/4~.$$ 
To estimate the product $\prod _{m=1}^{\infty}|1+x^2q^m|$ 
we define 
$l\in \mathbb{N}\cup 0$ by the condition $|q|^{-l}\leq |x|^2<|q|^{-l-1}$, 
$l\geq 6$. 
Thus for $m\geq l+2$ one has 

$$|x^2q^m|<|q^{m-l-1}|~~{\rm and}~~ 
\prod _{m=l+2}^{\infty}|1+x^2q^m|\geq \prod _{m=1}^{\infty}(1-|q|^m)=r~,
~~{\rm hence}$$ 
 
$$\prod _{m=1}^{l-1}|1+x^2q^m|\geq \prod _{m=1}^{l-1}(|q^{m-l}|-1)=|q|^{-l(l-1)/2}
\prod _{m=1}^{l-1}(1-|q|^{l-m})\geq |q|^{-l(l-1)/2}r~.$$ 
There remains to consider the product 

$$|1+x^{-2}q^l|\, |1+x^{-2}q^{l+1}|=(1/|x|^4)\, |x^2+q^l|\, |x^2+q^{l+1}|~.$$
The second and the third factor are not less than $\delta$ for 
$x^2\in (\mathbb{C}\backslash 0)\backslash \Sigma$. Hence the whole 
product is not less than $|q|^{2l+2}\delta ^2$. Thus 

$$|\Theta (q,x^2)|\geq r^3|q|^{-l(l-1)/2+2l+2}\delta ^2~.$$ 
The exponent of $|q|$ is negative for $l\geq 6$. For $h>0$ set 
$B_h:=\{ x\in \mathbb{C}\, |\, |x|^2\geq h\}$. Hence 
$|\Theta (q,x^2)|\geq r^3\delta ^2$ for $x\in B_{|q|^{-6}}\backslash \Sigma$. 
On the other hand $\max _{B_h}|\Xi (q,x^2)|$ 
tends to $0$ as $h\rightarrow \infty$. 
Hence for $h$ sufficiently large the equation $\theta (q,x^2)=0$, i.e. 
$\Theta ^*(q,x^2)=-\Xi (q,x^2)$, has no solution in 
$B_{|q|^{-6}}\backslash \Sigma$. Part (1) is proved.

%The proof of part (2) is analogous. The quantities $|\Xi |$ and $|\Xi '|$ are 
%uniformly small in the set $\Sigma$ for $|q|$ sufficiently small because 
%the moduli of the zeros $\mu _k$ of $\Theta ^*$ are not smaller than $1/|q|$. 
%Therefore there is a zero of $\theta$ close to $-1/q^k$ 
%for any $k\in \mathbb{N}$ (and not only for $k$ sufficiently large) 
%if $|q|$ is sufficiently small. 

To prove part (2) consider the function $h(q,x):=\prod (1+x/x_k)$, where 
$-x_k$ are the zeros of $\theta$ counted with multiplicity 
their moduli forming a non-decreasing sequence. (The convergence of 
the infinite product follows from part (1) -- for large values of $k$ 
the numbers $x_k$ are approximated by a geometric progression with 
ratio $1/q$.) Hence one has 
$\theta (q,x)=h(q,x)\Phi (q,x)$, where for each fixed 
$q$ the function $\Phi$ is an entire function without zeros. By Theorem~3 
of Chapter~I, Section~3 of \cite{Le} one has $\Phi =e^{\varphi}$, where for 
each fixed $q$, $\varphi$ is an entire function.

The order of the product of two entire functions of different (of equal) 
orders is the 
greater of the two orders 
(is their common order, see Theorem~12 of Chapter~I, Section~9 of \cite{Le}). 
The order of $\theta$ is defined by the formula (see 
Theorem~2 of Chapter~I, Section~2 of \cite{Le})

\begin{equation}\label{formula}
\overline{\lim _{k\rightarrow \infty}}(k\ln k/\ln (1/|q^{k(k+1)/2}|))=0~.
\end{equation}
Hence $\Phi$ and $h$ must be (for each fixed $q$) 
both of order $0$. This implies that $\varphi$ 
must be a constant. As $\theta (q,0)=1=h(q,0)$, one must have 
$\Phi =1$ (i.e. $\varphi =0$). This proves part (2) of the theorem.

Now we prove part (3).  
It is shown in \cite{Ko2} that for $q\in (\tilde{q}_j,\tilde{q}_{j+1}]$ 
the function $\theta (q,.)$ has $j$ conjugate pairs of roots 
counted with multiplicity. Hence for these values of $q$ one can set 
$\theta =P(q,x)\psi (q,x)$, where for each fixed $q$, 
$P(q,.)$ is a degree $2j$ polynomial in $x$ without real roots and 
the zeros of $\psi (q,.)$ (resp. of $P(q,.)$) 
are the real (resp. the complex) zeros of $\theta (q,.)$ counted with 
multiplicity. We assume that $P(q,0)=1$. Recall that 
the real zeros of $\theta$ are all negative.

%To prove part (2) of the theorem there remains to show that 
%$\psi (q,.)\in \cal{LP-I}$. 
Consider the infinite product 
$g(q,x):=\prod (1+x/\xi _k)$, where $-\xi _k$ are the real zeros of 
$\theta (q,.)$ 
given in the decreasing order. Theorem~4 in \cite{Ko2} implies that 
$\lim _{k\rightarrow \infty}\xi _kq^k=1$. Hence there exists $D>0$ such that 
$\xi _k\geq Dq^{-k}$. Set $g:=g_0+g_1x+g_2x^2+\cdots$. One has $g_k>0$ and  

\begin{equation}\label{gk}
g_k=\sum _{1\leq j_1<j_2<\cdots <j_k}1/\xi _{j_1}\xi _{j_2}\cdots \xi _{j_k}\leq D^{-k}S_k~,
\end{equation}
where $S_k:=\sum _{1\leq j_1<j_2<\cdots <j_k}q^{j_1+j_2+\cdots +j_k}$. It is clear that 
$S_1=q/(1-q)$ and 

$$S_k\leq \sum _{1\leq j_1<j_2<\cdots <j_{k-1}}q^{j_1+j_2+\cdots +j_{k-1}}
\sum _{j_k=k}^{\infty}q^{j_k}=S_{k-1}q^k/(1-q)~.$$
Hence by induction on $k$ one shows that $S_k\leq q^{k(k+1)/2}/(1-q)^k$ and  
$g_k\leq q^{k(k+1)/2}/(D(1-q))^k$. 

Thus for each fixed $q\in (0,1)$ the product $g$ is an entire function of 
the Laguerre-P\'olya class $\cal{LP-I}$. Its order equals $0$. 
Indeed, using again Theorem~2 of Chapter~I, Section~2 of \cite{Le},  
the order of $g$ 
is defined by the formula  

\begin{equation}\label{formula1}
\overline{\lim _{k\rightarrow \infty}}(k\ln k/\ln (1/|g_k|))~.
\end{equation} 
For large values of $k$ one has 

$$1/|g_k|\geq (D(1-q))^kq^{-k(k+1)/2}~~{\rm hence}~~
\ln (1/|g_k|)\leq -k(k+1)(\ln q)/2+k(\ln D+\ln (1-q))$$
and the above limit is $0$. 

The function $\psi$ equals $gg_1$, where $g_1$ 
is an entire function without zeros. As in the proof of part (1) one shows 
that $g_1=e^{\tilde{g}}$, where for each fixed $q$,  
$\tilde{g}$ is an entire function, and 
that the order of $g_1$ equals $0$. 
Hence $\tilde{g}$ is a constant and as $\psi (q,0)=g(q,0)=1$, one must have 
$\tilde{g}\equiv 0$, i.e. $g_1\equiv 1$.

Part (4) of the theorem results from part (3) -- for $k$ large enough the zero 
of $\theta$ which is close to $\mu _k=-q^{-k}$ is simple and only 
finitely-many of its zeros are not of this kind. 

To prove part (5) notice that for $q\in (-1,0)$ the numbers $\mu _k$ 
are real and change alternatively sign. As $\theta (q,x)=\theta (\bar{q},x)$, 
part of the zeros of $\theta$ are real and the rest form conjugate pairs. 
For large $k$ there is just one zero of $\theta$ close to $\mu _k$ 
hence this zero is real. The condition 
$\tilde{\xi}_k\tilde{\xi}_{k+1}<0$ follows from the alternative changing of 
sign of $\mu _k$. As in the proof of part (3) we show that 
$\theta =R\Lambda _1\Lambda _2$, where $R$ is a real polynomial with constant 
term 1 and without real zeros and $\Lambda _i$ (of the form 
$\prod _k(1+x/\tilde{\xi}_{i,k})$) are functions of the 
Laguerre-P\'olya class $\cal{LP-I}$ ($\tilde{\xi}_{1,k}>0$, 
$\tilde{\xi}_{2,k}<0$). Hence 
$\Lambda :=\Lambda _1\Lambda _2\in \cal{LP}$. For $k$ large enough the 
sequence $\{ |\tilde{\xi}_k|\}$ is monotone increasing because $\tilde{\xi}_k$ 
is close to $\mu _k$.    
Theorem \ref{basictm} is proved. 
%\end{proof}

\end{document}